\documentclass{elsart_mod}

\usepackage{amssymb}
\usepackage{amsmath} 
\usepackage{mathrsfs}
\usepackage[american]{babel}
\usepackage[PostScript=dvips]{diagrams}
\usepackage{cite}
\usepackage{url}
\usepackage{calligra}
\usepackage{rotating} 
\usepackage{bbold} 

\allowdisplaybreaks

\begin{document}

\newcommand{\zen}[1]{\mbox{\rm Z}(#1)}
\newcommand{\ZZ}{\mbox{$\mathbb{Z}$}}
\newcommand{\sZZ}{\mbox{$\scriptstyle\mathbb{Z}$}}
\newcommand{\CC}{\mbox{$\mathbb{C}$}}
\newcommand{\sCC}{\mbox{$\scriptstyle\mathbb{C}$}}
\newcommand{\FF}{\mbox{$\mathbb{F}$}}
\newcommand{\NN}{\mbox{$\mathbb{N}$}}
\newcommand{\QQ}{\mbox{$\mathbb{Q}$}}
\newcommand{\cls}[1]{\mathscr{C}_{#1}}
\newcommand{\chR}[2]{\!\!\text{\footnotesize\textcalligra{Ch}}\,\,_{#1}(#2)}
\newcommand{\clR}[2]{\!\!\text{\footnotesize\textcalligra{Cl}}\,\,_{#1}(#2)}
\newcommand{\bone}{\mathbb{1}}

\begin{frontmatter}


\title{On isomorphisms between centers of integral group rings of 
finite groups}
\runtitle{Isomorphisms between centers of integral group rings}


\author{Martin Hertweck} 
\ead{hertweck@mathematik.uni-stuttgart.de}
\address{Universit\"at Stuttgart, Fachbereich Mathematik,
IGT, 70569 Stuttgart, Germany}

\begin{abstract}
For finite nilpotent groups $G$ and $G^{\prime}$, and a $G$-adapted ring
$S$ (the rational integers, for example), it is shown 
that any isomorphism between the centers of the group rings $SG$ and
$SG^{\prime}$ is monomial, i.e., maps class sums in $SG$ to 
class sums in $SG^{\prime}$ up to multiplication with roots of unity.
As a consequence, $G$ and $G^{\prime}$ have identical character tables
if and only if the centers of their integral group rings $\ZZ G$ and 
$\ZZ G^{\prime}$ are isomorphic.
In the course of the proof, a new proof of the class sum correspondence 
is given.
\end{abstract}


\end{frontmatter}


\section{Introduction}\label{Sec:Intro}

It has been asked whether an automorphism of the center of the integral
group ring of a finite group necessarily induces a monomial permutation
on the set of the class sums
(listed as Problem~14.2 in the Kourovka Notebook \cite{Kouro:02}
and Problem~41 in \cite{Seh:93}, both times attributed to S.~D.~Berman),
but seemingly no progress was made towards a solution, except that 
it is annotated in \cite{Seh:93} that A.~A.~Bovdi has answered it
affirmatively for nilpotent groups of class at most three.
In this paper, it is finally dealt with the case of nilpotent groups.

Actually, we are treating an obvious generalization of the original question.
Suppose $G$ is a finite group and $S$ is an integral domain of 
characteristic zero. Let $\clR{S}{G}$ be the center of $SG$, the group ring 
of $G$ over $S$. Throughout, we shall be concerned with a $G$-adapted
coefficient ring $S$ (precise definitions are given in the following sections).
Then, given another group $G^{\prime}$, we ask whether an isomorphism
$\clR{S}{G}\cong\clR{S}{G^{\prime}}$ of $S$-algebras (if existing) is
necessarily monomial, i.e., maps class sums in $SG$ to 
class sums in $SG^{\prime}$ up to multiplication with roots of unity
(which, in any case, can be avoided by ``normalization,'' see
Lemma~\ref{mon}). For nilpotent $G$, this is answered in the affirmative 
(Theorem~\ref{nilthm}).
In fact, nilpotent groups constitute a special case as the familiar
Berman--Higman result assures that the central units of finite order
in $SG$ are, up to roots of unity, just the central elements in $G$.
This suggests that an approach might exist which proceeds inductively 
along the upper central series of the groups. This idea is realized,
resulting in an elementary, character free proof.

Dealing with isomorphisms---rather than only with automorphisms---%
does not cause additional difficulties and is motivated by the following 
question: Does an isomorphism
$\clR{\sZZ}{G}\cong\clR{\sZZ}{G^{\prime}}$ results in an isomorphism
between the character tables of $G$ and $G^{\prime}$?
(The converse is known ever since Frobenius introduced group characters.)
It turns out that a monomial isomorphism of centers preserves the
character degrees, thus giving rise to an isomorphism of character tables. 
Conversely, a degree preserving isomorphism of centers is monomial.
This is the content of the class sum correspondence, for which a proof is
given in Section~\ref{Sec:CSC} which deviates from the known in so far as 
the case when $S$ is a ring of algebraic integers is treated in a simple 
way, using an elementary result of Kronecker which says that at least one 
of the conjugates of a nonzero algebraic integer which is not a root of 
unity must lie outside the unit circle
(Lemmas~\ref{eps1} and \ref{eps2}).
It is only in Section~\ref{Sec:CSC} that characters really show up, 
and one might get an impression about the problem in general.
Section~\ref{Sec:Generals} contains all the generalities needed for
the handling of nilpotent groups in Section~\ref{Sec:Nilpotent}.
We would like to emphasize that no internal characterization of class sums,
or even character degrees, in the center is known, so we do {\em not}
consider $\clR{S}{G}$ as an algebra in its own right.

It seems appropriate to include a few remarks on character rings.
Let $\chR{\!S}{G}$ be the ring of $S$-linear combinations of the irreducible
characters of $G$. Let $R$ be a ring of algebraic integers.
Weidman \cite{Weid:65} and Saksonov \cite{Sak:66b} proved independently 
that if $\chR{\!R}{G}\cong\chR{\!R}{G^{\prime}}$,
then the character tables of $G$ and $G^{\prime}$ are the same.
This can be viewed as a consequence of an internal characterization
of the ordinary inner product on the character ring (see also the 
presentation in \cite{Ban:67}). Actually, they showed that any 
isomorphism of character rings is monomial (in the supplement \cite{Yam:96}
it is shown how ``normalization'' is to be understood). 

Duality between the rings $\clR{S}{G}$ and $\chR{\!S}{G}$ 
is definitive only for abelian groups. For a subfield $k$ of the complex
numbers, the structure of $\clR{k}{G}$ and $\chR{\!k}{G}$ has been
described in \cite{Tho:70}, where it is shown that if $p$ is an odd
prime and $G$ is a $p$-group, then $\clR{k}{G}\cong\chR{\!k}{G}$, 
with the assumption $p\neq 2$ being necessary (cf.\ also \cite{Ber:53}).
In \cite{Sak:79}, $p$-adic class algebras have been compared with 
$p$-adic character rings. If $S$ is $G$-adapted, then  
$\clR{S}{G}\cong\chR{\!S}{G^{\prime}}$ if and only if $G$ and 
$G^{\prime}$ are isomorphic abelian groups.

\section{Generalities}\label{Sec:Generals}

Let $G$ and $G^{\prime}$ be finite groups, and let $R$ be a ring
of algebraic integers in the field of complex numbers $\CC$. 
A {\em class sum} in $RG$ is, for a group element $g$ of $G$, the
sum of its $G$-conjugates in $RG$. Note that the class sums in $RG$
form an $R$-basis of $\clR{R}{G}$.

Suppose there exists an $R$-algebra isomorphism
$\varphi\colon\clR{R}{G}\rightarrow\clR{R}{G^{\prime}}$. In this section
we shall derive some basic facts about it.
Of course, $\varphi$ extends to a $\CC$-algebra isomorphism
$\varphi\colon\clR{\sCC}{G}\rightarrow\clR{\sCC}{G^{\prime}}$, and
as such, $\varphi$ maps primitive idempotents (corresponding to 
irreducible characters) to primitive idempotents and 
is completely determined by this operation.

We let $g_{1},\dotsc,g_{h}$ with $g_{1}=1$ be representatives of the classes
of $G$ and write $\mathscr{C}_{i}$ for the class sum of $g_{i}$ in $RG$.
We let $\varepsilon\colon RG\rightarrow R$ be the augmentation homomorphism
and $\varepsilon_{1}\colon RG\rightarrow R$ be the usual trace map.
These maps are defined, using the group basis $G$, by
$\varepsilon(\sum_{g\in G}a_{g}g)=\sum_{g\in G}a_{g}$ and 
$\varepsilon_{1}(\sum_{g\in G}a_{g}g)=a_{1}$ (all coefficients $a_{i}$
in $R$). 
Throughout, we dispose of similar notation for $G^{\prime}$ using dashes,
i.e., $g_{1}^{\prime},\dotsc,g_{h^{\prime}}^{\prime}$ are representatives of 
the classes of $G^{\prime}$ (with $g_{1}^{\prime}=1$), and so forth.
We also write $\mathscr{C}_{g}$ for the class sum of a specific element 
$g$ of $G$ (accordingly, $\mathscr{C}_{g^{\prime}}$ for $g^{\prime}$ in 
$G^{\prime}$ is defined) and let $|\mathscr{C}_{g}|$ be its length.

\begin{lem}\label{lin}
The groups $G$ and $G^{\prime}$ have the same order and the same number
of conjugacy classes. The isomorphism $\varphi$ maps the set of linear 
characters of $G$ onto the set of linear characters of $G^{\prime}$.
\end{lem}
\begin{pf}
Since $\clR{R}{G}$ and $\clR{R}{G^{\prime}}$ have the class sums as $R$-bases,
the groups $G$ and $G^{\prime}$ have the same number of conjugacy classes.
Without lost of generality we can assume that $|G|\geq |G^{\prime}|$.
Let $I$ be the ring of all algebraic integers in $\CC$. 
Let $\lambda$ be a linear character of $G$ and $\chi^{\prime}$
its image under $\varphi$. If $e_{\lambda}$ denotes the block idempotent 
in $\CC G$ corresponding to $\lambda$, then
\[ \bigg(\frac{|G^{\prime}|}{\chi^{\prime}(1)}e_{\lambda}\bigg)\varphi
= \sum_{i=1}^{h}\chi^{\prime}({g_{i}^{\prime}}^{-1})\mathscr{C}_{i}^{\prime}
\in \clR{I}{G^{\prime}}, \]
so $\frac{|G^{\prime}|}{\chi^{\prime}(1)}e_{\lambda}\in\clR{I}{G}$
and $\frac{|G^{\prime}|}{\chi^{\prime}(1)}\big/ |G|=
\varepsilon_{1}\big(\frac{|G^{\prime}|}{\chi^{\prime}(1)}e_{\lambda}\big)
\in I$ (here we used that $\lambda(1)=1$). Thus
$\frac{|G^{\prime}|}{\chi^{\prime}(1)}\big/ |G|$ is a natural integer,
and by our assumption on the orders of $G$ and $G^{\prime}$, this is only
possible if $\chi^{\prime}(1)=1$ and $|G|=|G^{\prime}|$.
\qed
\end{pf}

We shall say that $\varphi$ is {\em monomial} if there is a permutation
$\pi$ on $\{1,2,\dotsc,h\}$ and roots of unity 
$\xi_{1},\dotsc,\xi_{h}$ such that 
$\mathscr{C}_{i}\varphi=\xi_{i\pi}^{-1}\mathscr{C}_{i\pi}^{\prime}$ for
$1\leq i\leq h$. Note that then the assignment 
$\lambda^{\prime}(g_{i}^{\prime})=\xi_{i}$ defines
a linear character $\lambda^{\prime}$ of $G^{\prime}$ since the idempotent
$\frac{1}{|G|}\sum_{i=1}^{h}\mathscr{C}_{i}$ corresponding to the
principal character is mapped under $\varphi$ to 
$\frac{1}{|G^{\prime}|}\sum_{i=1}^{h}\xi_{i}^{-1}\mathscr{C}_{i}^{\prime}$.
We shall say that $\varphi$ is {\em normalized} if the following 
diagram is commutative:
\newcommand{\tina}%
{\begin{turn}{43}\raisebox{-2pt}{$\varepsilon$}\end{turn}}%
\newcommand{\tinb}%
{\begin{turn}{-43}\raisebox{5pt}{$\varepsilon^{\prime}$}\end{turn}}%
\begin{diagram}[width=2em,height=2em]
\clR{R}{G} & & \rTo^{\varphi} & & \clR{R}{G^{\prime}} \\
& \rdTo_{\tina} & & \ldTo_{\tinb} & \\
& & R & &
\end{diagram}
We write $\bone$ and $\bone^{\prime}$ for the principal characters
of $G$ and $G^{\prime}$ respectively. Note that 
$\varphi$ is normalized if and only if $\varphi$ sends
$\bone$ to $\bone^{\prime}$.
If $\varphi$ is monomial and normalized, it maps class sums to class sums.

We can turn our attention to normalized isomorphisms:

\begin{lem}\label{mon}
Let $\lambda^{\prime}$ be the linear character of $G^{\prime}$
to which $\bone$ is mapped under $\varphi$ (cf.\ Lemma~\ref{lin}).
Then a monomial $R$-algebra automorphism $\alpha$ of $\clR{R}{G^{\prime}}$
is defined by the assignment 
$\mathscr{C}_{i}^{\prime}\alpha=\lambda^{\prime}(g_{i}^{\prime})
\mathscr{C}_{i}^{\prime}$, and $\varphi\alpha$ is a normalized isomorphism.
\end{lem}
\begin{pf}
Note that $\varphi$ maps $\sum_{i}\mathscr{C}_{i}$ to 
$\sum_{i}\lambda^{\prime}({g_{i}^{\prime}}^{-1})\mathscr{C}_{i}^{\prime}$,
so $\alpha$ is clearly an $R$-linear map such that 
$(\sum_{i}\mathscr{C}_{i})\varphi\alpha=\sum_{i}\mathscr{C}_{i}^{\prime}$,
i.e., $\varphi\alpha$ sends $\bone$ to $\bone^{\prime}$.
It remains to show that $\alpha$ is multiplicative. Take class sums
$\mathscr{C}_{i}^{\prime}$ and $\mathscr{C}_{j}^{\prime}$ and write
$\mathscr{C}_{i}^{\prime}\mathscr{C}_{j}^{\prime} 
=\sum_{k}c_{ijk}^{\prime}\mathscr{C}_{k}^{\prime}$ with integers
$c_{ijk}^{\prime}$. Then 
\[ (\mathscr{C}_{i}^{\prime}\mathscr{C}_{j}^{\prime})\alpha=\sum_{k}
c_{ijk}^{\prime}\lambda^{\prime}(g_{k}^{\prime})\mathscr{C}_{k}^{\prime},
\quad (\mathscr{C}_{i}^{\prime})\alpha(\mathscr{C}_{j}^{\prime})\alpha=
\sum_{k}c_{ijk}^{\prime}\lambda^{\prime}(g_{i}^{\prime})
\lambda^{\prime}(g_{j}^{\prime})\mathscr{C}_{k}^{\prime}, \]
so we have to check that $\lambda^{\prime}(g_{k}^{\prime})=
\lambda^{\prime}(g_{i}^{\prime})\lambda^{\prime}(g_{j}^{\prime})$
whenever $c_{ijk}^{\prime}\neq 0$. But this is obvious by the definition
of the $c_{ijk}^{\prime}$ since $\lambda^{\prime}$ is a 
linear character, i.e., has the commutator subgroup of $G$
in its kernel. \qed
\end{pf}

We write $\widehat{N}$ for the sum in $RG$ of the elements of a 
subgroup $N$ of $G$.

\begin{lem}\label{norsgps}
Let $\varphi$ be normalized. Assume that there are 
normal subgroups $N$ and $N^{\prime}$ of $G$ and $G^{\prime}$, 
respectively, such that $\varphi$ maps
$\clR{R}{G}\cap RN$ onto $\clR{R}{G^{\prime}}\cap RN^{\prime}$.
Then $\widehat{N}\varphi=\widehat{N^{\prime}}$ and
there exists a commutative diagram of algebra homomorphisms
\begin{diagram}[width=4em,height=2em]
\clR{R}{G} & \rTo^{\varphi} & \clR{R}{G^{\prime}} \\ 
\dTo<{\pi} & & \dTo>{\pi^{\prime}} \\
\clR{R}{G/N} & \rTo^{\bar{\varphi}} & \clR{R}{G^{\prime}/N^{\prime}} 
\end{diagram}
where the vertical maps are the natural ones and the bottom map 
$\bar{\varphi}$ is again a normalized isomorphism.
\end{lem}
\begin{pf}
Let $e_{N}$ be the block idempotent of $\CC N$ corresponding to the
principal character. As a central idempotent of $\CC G$, it is
mapped by $\varphi$ to a central idempotent $f$ of $\CC G^{\prime}$.
By assumption, $f$ lies in $\CC N^{\prime}$; it is primitive in 
$\CC N^{\prime}$ as $e_{N}$ is primitive in $\CC N$. Since 
$\varphi$ is normalized, it follows that $f$ corresponds to the
principal character of $N^{\prime}$.
So $\widehat{N}\varphi=\frac{|N|}{|N^{\prime}|}\widehat{N^{\prime}}$,
showing that $\frac{|N|}{|N^{\prime}|}\in R\cap\QQ=\ZZ$.
Similarly, $\frac{|N^{\prime}|}{|N|}\in\ZZ$, so 
$|N|=|N^{\prime}|$ and $\widehat{N}\varphi=\widehat{N^{\prime}}$.

Set $\bar{G}=G/N$, and let $\pi\colon RG\rightarrow R\bar{G}$ be the
$R$-algebra homomorphism extending the natural homomorphism
$G\rightarrow\bar{G}$. Note that the kernel of $\pi$ is the annihilator
of $\widehat{N}$. Also considering the analogously defined map
$\pi^{\prime}$, it follows from $\widehat{N}\varphi=\widehat{N^{\prime}}$
that $\varphi$ induces an isomorphism
$\bar{\varphi}\colon\clR{R}{G}\pi\rightarrow\clR{R}{G^{\prime}}\pi^{\prime}$.
We proceed to show that $\bar{\varphi}$ can be extended to an isomorphism
$\clR{R}{\bar{G}}\rightarrow\clR{R}{\bar{G}^{\prime}}$.
Let $g\in G$, and let $\mathscr{C}_{\bar{g}}$ be the class sum of
$\bar{g}$ in $R\bar{G}$. Since $\pi$ maps conjugacy classes of $G$ onto
conjugacy classes of $\bar{G}$, we have 
$\mathscr{C}_{g}\pi=m\mathscr{C}_{\bar{g}}$ for some $m\in\NN$.
So $\clR{R}{G}\pi$ contains a $\CC$-basis of $\clR{\sCC}{\bar{G}}$, and
$\bar{\varphi}$ uniquely extends to an isomorphism
$\bar{\varphi}\colon\clR{\sCC}{\bar{G}}\rightarrow
\clR{\sCC}{\bar{G}^{\prime}}$.
With $\varphi$ also $\bar{\varphi}$ is normalized. It remains to show that 
$\mathscr{C}_{\bar{g}}\bar{\varphi}\in R\bar{G}^{\prime}$. Since
$\mathscr{C}_{\bar{g}}\bar{\varphi}=\frac{1}{m}
(\mathscr{C}_{g}\varphi\pi^{\prime})$, this is equivalent to
$\mathscr{C}_{g}\varphi\pi^{\prime}\in mR\bar{G}^{\prime}$ or
$(\mathscr{C}_{g}\widehat{N})\varphi\pi^{\prime}\in m|N|R\bar{G}^{\prime}$.
Since $\mathscr{C}_{g}\widehat{N}\in m\ZZ G$ we have 
$(\mathscr{C}_{g}\widehat{N})\varphi=(\mathscr{C}_{g}\varphi)
\widehat{N^{\prime}}\in mRG\cap (RG^{\prime})\widehat{N^{\prime}}$, so 
$(\mathscr{C}_{g}\varphi)\widehat{N^{\prime}}=mx\widehat{N^{\prime}}$
for some $x\in RG^{\prime}$, and it follows
$(\mathscr{C}_{g}\widehat{N})\varphi\pi^{\prime}=
((\mathscr{C}_{g}\varphi)\widehat{N^{\prime}})\pi^{\prime}=
(mx\widehat{N^{\prime}})\pi^{\prime}=m|N^{\prime}|(x\pi^{\prime})$
as desired.
\qed
\end{pf}

We define an anti-automorphism $\circ$ on $\CC G$ in the usual way
by $(\sum_{g\in G}a_{g}g)^{\circ}=\sum_{g\in G}\bar{a}_{g}g^{-1}$
where $\bar{a}_{g}$ denotes the complex conjugate of the number $a_{g}$.
Note that $\circ$ fixes each central primitive idempotent in $\CC G$.
This shows that $\varphi$ commutes with these anti-automorphisms
in the sense that $(x^{\circ})\varphi=(x\varphi)^{\circ^{\prime}}$
for all $x\in\CC G$. Also note that 
$\varepsilon_{1}(\mathscr{C}_{i}^{}\mathscr{C}_{i}^{\circ})=
|\mathscr{C}_{i}|$ for any index $i$ and 
$\varepsilon_{1}(\mathscr{C}_{i}^{}\mathscr{C}_{j}^{\circ})=0$
for distinct indices $i,j$.

\begin{lem}\label{eps1}
Suppose that $\varepsilon_{1}(\mathscr{C}_{g}^{}\mathscr{C}_{g}^{\circ})=
\varepsilon_{1}^{\prime}((\mathscr{C}_{g}^{}\mathscr{C}_{g}^{\circ})\varphi)$
for some $g\in G$. Write 
$\mathscr{C}_{g}\varphi=\sum_{i=1}^{h}a_{i}^{}\mathscr{C}_{i}^{\prime}$
with all $a_{i}$ in $R$ and suppose further that 
$a_{i_{0}}\neq 0$ for some index $i_{0}$ with 
$|\mathscr{C}_{i_{0}}^{\prime}|\geq |\mathscr{C}_{g}|$. Then 
$a_{i_{0}}$ is a root of unity and 
$\mathscr{C}_{g}\varphi=a_{i_{0}}^{}\mathscr{C}_{i_{0}}$.
\end{lem}
\begin{pf}
We have $\varepsilon_{1}(\mathscr{C}_{g}^{}\mathscr{C}_{g}^{\circ})=
|\mathscr{C}_{g}|$ and $\varepsilon_{1}^{\prime}((\mathscr{C}_{g}^{}
\mathscr{C}_{g}^{\circ})\varphi)=
\sum_{i=1}^{h}|a_{i}^{}|^{2}|\mathscr{C}_{i}^{\prime}|$.
By an elementary result due to Kronecker \cite{Kro:1857},
either $a_{i_{0}}$ is a root of unity or some algebraic conjugate of
$a_{i_{0}}$ has absolute value strictly greater than $1$.
From the assumptions, it follows that $a_{i_{0}}$ is a root of unity,
and that all other coefficients $a_{i}$ vanish. \qed
\end{pf}

The immediate consequence is:

\begin{lem}\label{eps2}
If $\varepsilon_{1}(z)=\varepsilon_{1}^{\prime}(z\varphi)$ for all
$z\in\clR{R}{G}$ then $\varphi$ is monomial.
\end{lem}
\begin{pf}
By Lemma~\ref{mon}, we can assume that $\varphi$ is normalized.
For $n\geq 0$, let $T_{n}$ be the set of class sums of elements of $G$ of 
length $n$ (so $T_{0}=\emptyset$), and let $R[T_{<n}]$ be the $R$-span
of $T_{0},\dotsc,T_{n-1}$.
We shall prove by induction on $n$ that $\varphi$ maps $T_{n}$ onto 
$T_{n}^{\prime}$. For $n=0$ this is an empty statement, so let $n\geq 1$.
Suppose that $T_{n}\neq\emptyset$, and let $g\in G$ with 
$\mathscr{C}_{g}\in T_{n}$. By the induction hypothesis, $\varphi$ maps
$R[T_{<n}]$ bijectively onto $R[T_{<n}^{\prime}]$. Thus $g$ satisfies the
hypotheses of Lemma~\ref{eps1}, from which we conclude that 
$T_{n}\varphi\subseteq T_{n}^{\prime}$. By symmetry, 
$T_{n}\varphi=T_{n}^{\prime}$, and we are done.
\end{pf}

\section{The class sum correspondence}\label{Sec:CSC}

For a group $G$, an integral domain $S$ of characteristic zero is called 
{\em $G$-adapted} if no prime divisor of the order of $G$ is invertible in 
$S$. In this section, we keep previous notation but instead of $R$ we take 
a $G$-adapted ring $S$ into consideration. 
It is only at first sight that this is a more general assumption, and we 
shall derive the class sum correspondence from results of the previous
section. Let $\psi\colon\clR{S}{G}\rightarrow\clR{S}{G^{\prime}}$ be an
$S$-algebra isomorphism.

First, we give the explicit formula of the matrix $A=(a_{ij})$ which 
describes $\psi$ with respect to the bases formed by the class sums:
\[ \mathscr{C}_{i}^{}\psi=\sum_{j=1}^{h}a_{ij}\mathscr{C}_{j}^{\prime}
\quad (1\leq i\leq h). \]
Let $\chi_{1},\dotsc,\chi_{h}$ with $\chi_{1}=\bone$ be the irreducible
characters of $G$, and let $e_{l}$ be the block idempotent
corresponding to $\chi_{l}$ (so $e_{l}=\frac{\chi_{l}(1)}{|G|}\sum_{i=1}^{h}
\chi_{l}(g_{i}^{-1})\mathscr{C}_{i}$). We have 
$e_{l}^{}\psi=e_{l\sigma}^{\prime}$ for a permutation $\sigma$
on $\{1,2,\dotsc,h\}$. 
Writing both sides as linear combinations of class sums, we find that
\begin{equation}\label{eqA2}
(\chi_{l}(g_{1}),\dotsc,\chi_{l}(g_{h}))A=
\frac{\chi_{l\sigma}^{\prime}(1)}{\chi_{l}(1)}
(\chi_{l\sigma}^{\prime}(g_{1}^{\prime}),\dotsc,
\chi_{l\sigma}^{\prime}(g_{h}^{\prime})).
\end{equation}
To put it in matrix form, let $M$ be the monomial matrix whose
$(l,l\sigma)$ entry is $\chi_{l\sigma}^{\prime}(1)/\chi_{l}(1)$ (and all
other entries are zero), let $X$ be the character table of $G$ 
regarding the fixed orders on classes and characters 
(by convention, $X^{\prime}$ is the character table of $G^{\prime}$).
Then $XA=MX^{\prime}$. Solving for $A$ using the first orthogonality
relation yields
\begin{equation}\label{eqA3}
a_{ij} = \frac{1}{|G|} \sum_{l=1}^{h}
\frac{|\mathscr{C}_{i}|\chi_{l}(g_{i}^{-1})}{\chi_{l}(1)}
\chi_{l\sigma}^{\prime}(1)\chi_{l\sigma}^{\prime}(g_{j}^{\prime}).
\end{equation}

We can assume that the quotient field of $S$ is embedded in a field 
containing $\CC$, so that $a_{ij}\in\QQ(\zeta)$ where $\zeta$ is a 
complex primitive $|G|$-th root of unity. We remark that
the Galois group of $\QQ(\zeta)$ over $\QQ$ acts on the entries of $A$, for
if $\tau$ is a Galois automorphism, with $\zeta^{\tau}=\zeta^{n}$ for
$n\in\NN$ (coprime to $|G|$), then 
\[ a_{ij}^{\tau} = \frac{1}{|G|} \sum_{l=1}^{h}
\frac{|\mathscr{C}_{g_{i}^{n}}|\chi_{l}(g_{i}^{-n})}{\chi_{l}(1)}
\chi_{l\sigma}^{\prime}(1)\chi_{l\sigma}^{\prime}({g_{j}^{\prime}}^{n}). \]

The following remark, due to Saksonov, shows that all the lemmas from the
previous section hold with $R$ replaced by the more general ring $S$.

\begin{rem}\label{Gada}{\normalfont
The Galois action mentioned before 
shows in particular that all algebraic conjugates of $a_{ij}$ lie in $S$,
and this implies that all $a_{ij}$ are algebraic integers,
by a lemma we can attribute to Saksonov \cite[p.~190]{Sak:71}
(Lemma~3.2.2 in \cite{Kar:89}).
Thus if $R$ is the ring generated over $\ZZ$ by the entries of $A$ and
$A^{-1}$, then $R$ is a ring of algebraic integers, and $\psi$ restricts
to an isomorphism $\varphi\colon\clR{R}{G}\rightarrow\clR{R}{G^{\prime}}$.
}\end{rem}

We are now in a position to give a quick proof of the class sum
correspondence, which was treated by Berman, Glauberman, Passman and Saksonov
(names that should be mentioned at least, cf.\ \cite[(1.1)]{RoSc:87}).
Glauberman and Passman \cite[Theorem~C]{Pas:65} treated the case when
the coefficients are algebraic integers while Saksonov's version 
\cite{Sak:71} is for $G$-adapted coefficient rings.
Short proofs treating the case when $\ZZ$ serves as coefficient ring 
can be found in \cite[(3.17)]{Isa:94}, \cite[Chapter~14, Lemma~2.3]{Pas:77},
\cite{San:86}, \cite[(36.5)]{Seh:93}.

We say that $\psi$  {\em preserves the character degrees} if 
$\chi_{l}(1)=\chi_{l\sigma}^{\prime}(1)$ for $1\leq l\leq h$, or, 
equivalently, if $\psi$ can be extended to an isomorphism 
between the complex group rings $\CC G$ and $\CC G^{\prime}$.

We state the class sum correspondence following Saksonov
\cite{Sak:71} (as reported in \cite[Theorem~3.5.8]{Kar:89}).

\begin{thm}[Class sum correspondence]\label{CSC}
Let $\psi\colon\clR{S}{G}\rightarrow\clR{S}{G^{\prime}}$ be an $S$-algebra
isomorphism. Then $\psi$ is monomial if and only if
it preserves the character degrees.
\end{thm}
\begin{pf}
By Remark~\ref{Gada}, we can assume that $S=R$ and $\psi=\varphi$ as before.
Suppose that $\varphi$ is monomial. Then the first column of $A$ is the
transposed of $(1,0,\dotsc,0)$, so comparing the first entries in
\eqref{eqA2} shows that $\varphi$ preserves the character degrees.
Conversely, suppose that $\varphi$ preserves the character degrees.
Using the second orthogonality relation, \eqref{eqA3} shows that 
$a_{i1}=0$ for $i>1$, so $\varphi$ satisfies the assumption of 
Lemma~\ref{eps2} and is therefore monomial.  \qed
\end{pf}

\begin{rem}\label{chdeg}{\normalfont
We do not know whether, in general, $\psi$ preserves the character degrees.
What is immediate from \eqref{eqA2} is that the
$\chi_{l\sigma}^{\prime}(1)^{2}/\chi_{l}(1)$ are integers, so that
$\chi_{l\sigma}^{\prime}(1)$ and $\chi_{l}(1)$ have the same prime divisors
(consider also $\psi^{-1}$).
The question is whether heights of irreducible characters are preserved.
In \cite[\S~0]{Rob:01} it is pointed out that it does not seem
immediately obvious that the $p$-defects (or heights) of 
irreducible characters not of height zero can be determined from 
the knowledge of the isomorphism type of the center of the group ring
of $G$ over a $p$-adic ring alone.

We remark that the bilinear form 
$(x,y)=\varepsilon_{1}^{\prime}((x\psi)(y\psi)^{\circ^{\prime}})$
on $\clR{S}{G}$ is given with respect to the basis formed by the 
class sums by the matrix $AD^{\prime}A^{\ast}$, where
$D^{\prime}=\text{\rm diag}
(|\mathscr{C}_{1}^{\prime}|,\dotsc,|\mathscr{C}_{h}^{\prime}|)$
and $A^{\ast}$ is the hermitian transpose of $A$. The $(i,j)$ entry
is given by
\[ (AD^{\prime}A^{\ast})_{ij}=
\frac{1}{|G|} \sum_{l=1}^{h}\chi_{l\sigma}^{\prime}(1)^{2}
\frac{|\mathscr{C}_{i}|\chi_{l}(g_{i}^{-1})}{\chi_{l}(1)}
\frac{|\mathscr{C}_{j}|\chi_{l}(g_{j})}{\chi_{l}(1)}. \] %
}\end{rem}

\begin{rem}\label{idchtbl}{\normalfont
The groups $G$ and $G^{\prime}$ have identical character tables if
the isomorphism $\psi\colon\clR{S}{G}\rightarrow\clR{S}{G^{\prime}}$
is monomial and degree preserving. This is easily seen as follows.
Choose an ordering of the irreducible characters such that $\psi$ maps the
block idempotent $e_{i}$ belonging to the character $\chi_{i}$ to the block 
idempotent $e_{i}^{\prime}$ belonging to the character $\chi_{i}^{\prime}$.
By definition of monomial isomorphism, there is a linear character 
$\lambda^{\prime}$ of $G^{\prime}$ such that, after suitable ordering of the
classes, $\mathscr{C}_{j}\psi=\lambda^{\prime}(g_{j}^{\prime})
\mathscr{C}_{j}^{\prime}$,  and by Lemma~\ref{mon}, 
$|\mathscr{C}_{j}|=|\mathscr{C}_{j}^{\prime}|$. Since $\psi$ is
degree preserving, we can set 
$q_{ij}=|\mathscr{C}_{j}|/\chi_{i}(1)=
|\mathscr{C}_{j}^{\prime}|/\chi_{i}^{\prime}(1)$, and obtain
\[ q_{ij}\chi_{i}(g_{j})e_{i}^{\prime}=(q_{ij}\chi_{i}(g_{j})e_{i})\psi=
(\mathscr{C}_{j}e_{i})\psi=
\lambda^{\prime}(g_{j}^{\prime})\mathscr{C}_{j}^{\prime}e_{i}^{\prime}=
\lambda^{\prime}(g_{j}^{\prime})q_{ij}
\chi_{i}^{\prime}(g_{j}^{\prime})e_{i}^{\prime}. \] So
$\chi_{i}(g_{j})=(\lambda^{\prime}\otimes\chi_{i}^{\prime})(g_{j}^{\prime})$,
and the $\lambda^{\prime}\otimes\chi_{i}^{\prime}$ are, of course,
again the irreducible characters of $G^{\prime}$.
}\end{rem}

\section{Nilpotent groups}\label{Sec:Nilpotent}

We keep the notion introduced in the previous sections.
The following lemma is designed for application to nilpotent groups.

\begin{lem}\label{appl}
Keep notation and hypothesis from Lemma~\ref{norsgps}. 
Let $M$ be the normal subgroup of $G$ containing $N$ so that $M/N$ is the
center of $G/N$, and define the normal subgroup $M^{\prime}$ of $G^{\prime}$
analogously. Suppose that $\varphi$ maps the class sums of elements of $N$ 
onto the class sums of elements of $N^{\prime}$. Then $\varphi$ maps the 
class sums of elements of $M$ onto the class sums of elements of $M^{\prime}$.
\end{lem}
\begin{pf}
Set $\bar{G}=G/N$ and $\bar{G}^{\prime}=G^{\prime}/N^{\prime}$.
By the Berman--Higman result (from \cite{Ber:55} and \cite[p.~27]{Hig:40})
and Lemma~\ref{norsgps}, we know that $\bar{\varphi}$ maps the class sums of
elements of $\bar{M}$ (i.e., the central elements of $\bar{G}$)
onto the class sums of elements of $\bar{M}^{\prime}$
(i.e., the central elements of $\bar{G}^{\prime}$).

For $n\geq 0$, let $T_{n}$ be the set of class sums of elements of $M$ of 
length $n$ (so $T_{0}=\emptyset$), and define $T_{n}^{\prime}$ analogously.
We shall prove by induction on $n$ that $\varphi$ maps $T_{n}$ onto 
$T_{n}^{\prime}$. For $n=0$ this is an empty statement, so let $n\geq 1$.
Suppose that $T_{n}\neq\emptyset$, and let $m\in M$ with 
$\mathscr{C}_{m}\in T_{n}$.
By the Berman--Higman result, there is $m^{\prime}\in M^{\prime}$ with
$\mathscr{C}_{\bar{m}}\bar{\varphi}=\mathscr{C}_{\bar{m}^{\prime}}$, and
we can write
\begin{equation}\label{eq1}
\mathscr{C}_{m}\varphi=
\underset{\bar{m}_{i}^{\prime}=\bar{m}^{\prime}\text{ for all $i$}}
{\bigg(\underbrace{\sum_{i=1}^{k}r_{i}^{}\mathscr{C}_{m_{i}^{\prime}}}
\bigg)} +
\underset{\text{in the kernel of $\pi^{\prime}$}}{\bigg(
\underbrace{\sum_{j=1}^{l}s_{j}^{}\mathscr{C}_{h_{j}^{\prime}}}\bigg)} 
\end{equation}
with the class sums 
$\mathscr{C}_{m_{1}^{\prime}},\dotsc,\mathscr{C}_{m_{k}^{\prime}},
\mathscr{C}_{h_{1}^{\prime}},\dotsc,\mathscr{C}_{h_{l}^{\prime}}$
pairwise distinct. Since $\bar{m}\in\zen{\bar{G}}$ we have 
$\mathscr{C}_{m}^{}\mathscr{C}_{m}^{\circ}\in RN$ and so
$\varepsilon_{1}(\mathscr{C}_{m}^{}\mathscr{C}_{m}^{\circ})=
\varepsilon_{1}^{\prime}((\mathscr{C}_{m}^{}\mathscr{C}_{m}^{\circ})\varphi)$
by assumption, meaning that
\begin{equation}\label{eq2}
|\mathscr{C}_{m}|=
\bigg(\sum_{i=1}^{k}|r_{i}|^{2}|\mathscr{C}_{m_{i}^{\prime}}|\bigg) +
\bigg(\sum_{j=1}^{l}|s_{j}|^{2}|\mathscr{C}_{h_{j}^{\prime}}|\bigg).
\end{equation}
By Lemma~\ref{eps1} (application of Kronecker's result), if
$|\mathscr{C}_{m_{i}^{\prime}}|\geq |\mathscr{C}_{m}|$ for some index $i$
then $\mathscr{C}_{m}\varphi=\mathscr{C}_{m_{i}^{\prime}}$. So let us assume
that $|\mathscr{C}_{m_{i}^{\prime}}|<|\mathscr{C}_{m}|$ for all $i$.
We will reach a contradiction, showing that 
$T_{n}\varphi\subseteq T_{n}^{\prime}$. By symmetry, then also 
$T_{n}^{\prime}\varphi^{-1}\subseteq T_{n}$ and the proof will be complete.
By the induction hypothesis, there are $m_{i}\in M$ with 
$\mathscr{C}_{m_{i}}\varphi=\mathscr{C}_{m_{i}^{\prime}}$. Then 
$\bar{m}_{i}^{}\bar{\varphi}=\bar{m}_{i}^{\prime}=\bar{m}^{\prime}
=\bar{m}^{}\bar{\varphi}$ shows that $Nm=Nm_{i}$ for all $i$. Set
$\Delta=\mathscr{C}_{m}-\sum_{i=1}^{k}r_{i}\mathscr{C}_{m_{i}}$. Then
$\Delta\varphi=\sum_{j=1}^{l}s_{j}^{}\mathscr{C}_{h_{j}^{\prime}}$.
Since $\Delta\Delta^{\circ}\in RN$, we have
$\varepsilon_{1}(\Delta^{}\Delta^{\circ})=\varepsilon_{1}^{\prime}
((\Delta\varphi)(\Delta\varphi)^{\circ^{\prime}})$ which gives
\begin{equation}\label{eq3}
|\mathscr{C}_{m}|+
\bigg(\sum_{i=1}^{k}|r_{i}|^{2}|\mathscr{C}_{m_{i}}|\bigg) =
\bigg(\sum_{j=1}^{l}|s_{j}|^{2}|\mathscr{C}_{h_{j}^{\prime}}|\bigg).
\end{equation}
From \eqref{eq2} and \eqref{eq3} it follows that all $r_{i}$ are zero,
so that \eqref{eq1} gives the desired contradiction
$|\mathscr{C}_{m}|\bar{m}^{\prime}=\mathscr{C}_{m}\pi\bar{\varphi}=
\mathscr{C}_{m}\varphi\pi^{\prime}=0$.
\qed
\end{pf}

\begin{thm}\label{nilthm}
Let $G$ and $G^{\prime}$ be finite nilpotent groups, and let $S$ be a 
$G$-adapted ring. Then all isomorphisms between 
$\clR{S}{G}$ and $\clR{S}{G^{\prime}}$ as $S$-algebras (if there are any) 
are monomial and preserve the character degrees.
\end{thm}
\begin{pf}
By  Remark~\ref{Gada}, we can assume that $S=R$ and $\psi=\varphi$ as above.
By Lemma~\ref{mon} we can further assume that $\varphi$ is normalized,
and then we only need to show that $\varphi$ maps 
class sums to class sums, by Theorem~\ref{CSC}.

Set $Z_{n}=\text{\rm Z}_{n}(G)$, the $n$-th term of the upper central series
of $G$ (so $Z_{0}=1$ and $Z_{1}=\zen{G}$). Use similar notation for 
$G^{\prime}$. We prove by induction on $n$ that 
$\varphi$ maps the class sums of elements of $Z_{n}$ onto the class sums 
of elements of $Z_{n}^{\prime}$. For $n=0$, there is nothing to prove,
so we can let $n\geq 1$ when the statement follows from the induction 
hypothesis and Lemma~\ref{appl}, applied with $N=Z_{n-1}$ and $M=Z_{n}$
(and the corresponding normal subgroups of $G^{\prime}$). 
\qed
\end{pf}

By Remark~\ref{idchtbl}, we obtain as corollary:

\begin{cor}\label{c4}
Finite nilpotent groups $G$ and $G^{\prime}$ have identical character tables
if and only if $\clR{\sZZ}{G}\cong\clR{\sZZ}{G^{\prime}}$ as rings.
\end{cor}

\bibliographystyle{amsplain}

\providecommand{\bysame}{\leavevmode\hbox to3em{\hrulefill}\thinspace}

\bibliography{HeSo_BrPairs}

\end{document}